\documentclass{article}
\usepackage{fullpage}
\usepackage{amsmath,amsthm,amsfonts}

\def\Ext{\mathrm{Ext}}
\def\Hom{\mathrm{Hom}}
\def\mod{\mathrm{mod}}
\def\C{\mathcal{C}}
\def\Y{\mathcal{Y}}
\def\D{\mathcal{D}}
\def\ad{\mathrm{ad}}
\def\id{\mathrm{id}}
\def\wh{\widehat}
\def\ot{\otimes}

\newtheorem{teo}{Theorem}[section]
\newtheorem{defi}[teo]{Definition}
\newtheorem{lem}[teo]{Lemma}

\newtheorem{coro}[teo]{Corollary}

\author{Marco A. Farinati $^{1}$ - Andrea Solotar ${}^{1}$}

\title{G-structure on the cohomology of Hopf algebras}

\date{}
\begin{document}
\maketitle
\begin{footnotetext}[1]
{Dto. de Matem\'atica Facultad de Cs. Exactas y Naturales.
Universidad de Buenos Aires. Ciudad Universitaria Pab I. 1428 - 
Buenos Aires - Argentina. e-mail: asolotar@dm.uba.ar,  mfarinat@dm.uba.ar\\
Research partially supported by UBACYT X062 and
Fundaci\'on Antorchas (proyecto 14022-47).\\
Both authors are
research members of CONICET (Argentina).}
\end{footnotetext}

\begin{abstract}
We prove that $\Ext^{\bullet}_A(k,k)$ is a Gerstenhaber algebra, where
$A$ is a Hopf algebra. In case $A=D(H)$ is the Drinfeld double
of a finite dimensional
Hopf algebra $H$, our results implies the existence of a Gerstenhaber
bracket on $H^{\bullet}_{GS}(H,H)$. This fact was conjectured by
R. Taillefer in \cite{T3}. The method consists in identifying
$\Ext^{\bullet}_A(k,k)$ as a Gerstenhaber
subalgebra of  $H^{\bullet}(A,A)$ (the Hochschild cohomology of $A$).
\end{abstract}


\section*{Introduction}
The motivation of this paper is to prove that $H^{\bullet}_{GS}(H,H)$ 
has a structure of a G-algebra. We prove
this result
when $H$ is a finite dimensional 
Hopf algebra 
(see Theorem \ref{teo3} and Corollary \ref{coroimportante}). 
$H^{\bullet}_{GS}$ is the cohomology theory for Hopf algebras
defined by Gerstenhaber and Schack in \cite{GS1}. 

In order to obtain commutativity of the cup product 
we prove a general statement on $\Ext$ groups over 
Hopf algebras (without any finiteness assumption). When $H$ is finite
dimensional, the category of Hopf bimodules is isomorphic
to a module category, over an algebra $X$ (also finite dimensional)
defined by C. Cibils and M. Rosso (see \cite{CR}), and 
this category is also equivalent
to the category of Yetter-Drinfeld modules,
which is isomorphic to the category of modules over the Hopf algebra
$D(H)$ (the Drinfeld double of $H$). In \cite{T2}, R. Taillefer
has defined a natural cup product
in $H^{\bullet}_{GS}(H,H)=H^{\bullet}_b(H,H)$ (see \cite{GS2} 
for the definition of $H^{\bullet}_b$).  When $H$ is 
finite dimensional she proved that 
$H^{\bullet}_b(H,H)\cong\Ext^{\bullet}_{X}(H,H)$, and using this isomorphism
she showed that it is (graded) commutative. In a later work
\cite{T3} she extended the result of commutativity of the cup product
to arbitrary dimensional Hopf algebras and she conjectured
the existence (and a formula) of a Gerstenhaber bracket.

Our method for giving a Gerstenhaber bracket is the following:
under the equivalence
of categories ${}_X$-$\mod\cong {}_{D(H)}$-$\mod$, the object $H$ corresponds
to $H^{coH}=k$, so $\Ext^{\bullet}_{X}(H,H)\cong 
\Ext^{\bullet}_{D(H)}(k,k)$ (isomorphism
of graded algebras); after D. \c Stefan \cite{St} one knows that
$\Ext^{\bullet}_{D(H)}(k,k)\cong H^{\bullet}(D(H),k)$. 
In Theorem \ref{teo2}
we prove that, if $A$ is an arbitrary Hopf algebra, then
$H^{\bullet}(A,k)$ is isomorphic to a subalgebra of $H^{\bullet}(A,A)$ (in
particular it is graded commutative) and in Theorem \ref{teo3}
we prove that the image of 
$H^{\bullet}(A,k)$ in 
$H^{\bullet}(A,A)$ is stable under the brace operation, in particular
it is closed under the Gerstenhaber bracket of $H^(A,A)$. So,
the existence of the Gerstenhaber bracket on $H^{\bullet}_{GS}(H,H)$ follows,
at least in the finite dimensional case, taking $A=D(H)$. We don't
know if this bracket coincides with the formula
proposed in \cite{T3}.

We also provide a proof that the algebra 
$\Ext^{\bullet}_{\C}(k,k)$ is graded commutative when 
$\C$ is a braided monoidal category satisfying certain
homological hypothesis (see Theorem \ref{teo1}). This gives
an alternative proof of the commutativity result in the
arbitrary dimensional case taking $\C={}_H^H\Y\D$, the
Yetter-Drinfeld modules.

In this paper, the letter $A$ will denote a Hopf algebra over a 
field $k$.

%

\section{Cup products}

This section has two parts. First we prove a
generalization of the fact that the cup product
on $H^{\bullet}(G,k)$ is graded commutative. The general abstract setting
is that of a braided (abelian) category with enough injectives
satisfying a K\"unneth formula (see definitions below).
The other part will concern
the relation between self extensions of $k$ and Hochschild
cohomology of $A$ with coefficients in $k$. 

Let us recall the definition of a braided category:
\begin{defi}
The data $(\C,\ot,k,c)$ is called a {\bf braided} category 
with unit element $k$ if
\begin{enumerate}
\item $\C$ is an abelian category.
\item $-\ot -$ is a bifunctor, bilinear, associative, and
there are natural isomorphisms $k\ot X\cong X\cong X\ot k$
for all objects $X$ in $\C$.
\item For all pair of objects $X$ and $Y$, 
$c_{X,Y}:X\ot Y\to Y\ot X$ is a natural isomorphism. The isomorphisms
$c_{X,k}:X\ot k\cong k\ot X$ agrees with the isomorphism
of the unit axiom, and for all triple $X$, $Y$, $Z$
of objects in $\C$, the Yang-Baxter equation is satisfied:
\[
(\id_Z\ot c_{X,Y})\circ (c_{X,Z}\ot \id_Y)\circ(\id_X\ot c_{Y,Z})=
(\id_Y\ot c_{X,Z})\circ (c_{X,Y}\ot \id_Z)
\]
\end{enumerate}
If one doesn't have the data $c$, and axioms 1 and
2 are
satisfied,
we say that $(\C,\ot,k)$ is a {\bf monoidal} category.
\end{defi}

\begin{defi}
We will say that a monoidal category $(\C,\ot,k)$ satisfies
the {\bf K\"unneth formula} if and only if
there are natural isomorphisms
$H_*(X_*,d_X)\ot H_*(Y_*,d_Y)\cong H_*(X_*\ot Y_*,d_{X\ot Y})$
for all pair of complexes in $\C$.
\end{defi}

\begin{teo}\label{teo1}
Let $(\C, \ot,k,c)$ be a braided category with enough
injectives satisfying the K\"unneth
formula, then $\Ext^{\bullet}_{\C}(k,k)$ is graded commutative.
\end{teo}
\begin{proof} We proceed as in the proof that
$H^{\bullet}(G,k)$ is graded commutative (see for example
\cite{B}, page 51, Vol I).
The proof is based on two points: firstly a definition
of a cup product using $\ot$,
secondly a Lemma relating this construction and
the Yoneda product of extensions.

Let $0\to M\to X_p\to \dots X_1\to N\to 0$
and  $0\to M'\to X'_q\to \dots X'_1\to N'\to 0$ be
two extensions in $\C$. Then
 $N_*:=(0\to M\to X_p\to \dots X_1\to 0)$
and  $N'_*:=(0\to M'\to X'_q\to \dots X'_1\to 0)$ are two complexes,
quasi-isomorphic to $N$ and $N'$ respectively.
By the K\"unneth formula $N_*\ot N'_*$ is a complex quasi-isomorphic
to $N\ot N'$, so "completing" this complex with $N\ot N'$
(more precisely considering the mapping cone of the chain map
$N_*\ot N'_*\to N\ot N'$) one has an extension in $\C$,
beginning with $M\ot M'$ and ending with $N\ot N'$.

So, we have defined a cup product:
\[\Ext^p_{\C}(N,M)\times\Ext_{\C}^q(N', M')\to
\Ext_{\C}^{p+q}(N\ot N',M\ot M')\]
We will denote this product by a dot, and the
Yoneda product by $\smile$.
The Lemma relating this product and the Yoneda one
is the following:
\begin{lem}If $f\in\Ext^p_{\C}(M,N)$ and
$g\in\Ext^q_{\C}(M',N')$, then
\[f\smile g=(f\ot \id_{N'})\smile (\id_M\ot g)\]
\end{lem}
{\em Proof of the Lemma:} Interpreting the elements $f$ and
$g$ as extensions, it is clear how to define a 
morphism of complexes
$(f\ot \id_{N'})\smile (\id_M\ot g)\to f\smile g$,
and by the K\"unneth formula, it is a quasi-isomorphism.

In the particular case $M=M'=N=N '=k$, the Lemma implies
that $f.g=f\smile g$ for all $f$ and $g$ in $\Ext^{\bullet}_{\C}(k,k)$. 
Now the theorem is a consequence
of the isomorphism
$(X_*\ot Y_*,d_{X\ot Y})\cong 
(Y_*\ot X_*,d_{Y\ot X})$, valid for every pair of complexes
in $\C$,
defined by:
\[(-1)^{pq}c_{X,Y}:X_p\ot Y_q\to Y_q\ot X_p\] 
\end{proof}
\begin{teo}\label{teo2}
If $A$ is a Hopf algebra 
then 
$\Ext^{\bullet}_{\C}(k,k)\cong H^{\bullet}(A,k)$. Moreover
$H^{\bullet}(A,k)$ is isomorphic to  a 
subalgebra of
$H^{\bullet}(A,A)$.
\end{teo}
\begin{proof}  After D. \c Stefan \cite{St}, since $A$ is an $A$-Hopf
Galois extension of $k$, $H^{\bullet}(A, M)\cong \Ext^{\bullet}_A(k,M^{\ad})$
for all $A$-bimodule $M$.
In particular, $H^{\bullet}(A,k)=\Ext^{\bullet}_A(k,k)$.
But one can give, for this particular case, 
an explicit morphism at the complex level. In order
to do this, we will choose a particular resolution
of $k$ as left $A$-module.

Let $C_*(A,b')$ be the standard resolution
of $A$ as $A$-bimodule, namely
$C_n(A,b')=A\ot A^{\ot n}\ot A$ and
$b'(a_0\ot\dots \ot a_{n+1})=
\sum_{i=0}^n(-1)^{i}a_0\ot\dots \ot a_i.a_{i+1}\ot
\dots \ot a_{n+1}$ ($a_i\in A$). This resolution splits on the right,
so $(C_*(A)\ot_A k,b'\ot \id_k)$ is a resolution
of $A\ot_A k=k$ as left $A$-module.
Using this resolution, $\Ext^{\bullet}_A(k,k)$ is the homology
of the complex $(\Hom_A(C_*(A)\ot_Ak,k), (b'\ot_A\id_k)^{*})\cong
(\Hom(A^{\ot *},k),\partial)$. Under this isomorphism,
the differential
$\partial$ is given by
\[(\partial f)(a_1\ot\dots\ot a_n)=\epsilon(a_1)f(a_2\ot\dots \ot a_n)+\]
\[+
\sum_{i=1}^{n-1}(-1)^if(a_1\ot\dots\ot a_i.a_{i+1}\ot\dots\ot a_n)+
(-1)^nf(a_1\ot\dots\ot a_{n-1})\epsilon(a_n)\]
And this is precisely the formula of the differential
of the standard Hochschild complex
computing
$H^{\bullet}(A,k)$.

One can easily check that the cup product on $\Ext^{\bullet}_A(k,k)$
(which equals the Yoneda product in this case) corresponds
to the cup product on $H^{\bullet}(A,k)$, so this isomorphism is an 
algebra isomorphism.

Now we will give two multiplicative maps $H^{\bullet}(A,k)\to H^{\bullet}(A,A)$
and $H^{\bullet}(A,A)\to H^{\bullet}(A,k)$.
Consider the counit $\epsilon:A\to k$, it is an algebra map, so
the induced map $\epsilon_*:H^{\bullet}(A,A)\to H^{\bullet}(A,k)$ is
multiplicative. We will define a multiplicative
section of this map.

Let $f:A^{\ot p}\to k$ be a Hochschild cocycle, define
$F:A^{\ot p}\to A$ by the formula:
\[
F(a^1\ot\dots\ot a^p):=
a^1_1 \dots a^p_1.
f(a^1_2\ot\dots\ot a^p_2)
\]
Where we have used the Sweedler-type notation with
summation symbol omitted: $a^i_1\ot a^i_2=\Delta(a^i)$,
for $a^i\in A$.

Let us check that $F$ is a Hochschild cocycle with values
in $A$.
\[
\partial(F)(a^0\ot\dots \ot a^p)=
a^0F(a^1\ot\dots \ot a^p)+\]
\[+
\sum_{i=0}^{p-1}(-1)^{i+1}F(a^0\ot\dots\ot a^i.a^{i+1}\ot\dots  \ot a^p)
+(-1)^{p+1}F(a^0\ot\dots \ot a^{p-1})a^p= \]
\[=a^0.a^1_1\dots a^p_1.f(a^1_2\ot\dots \ot a^p_2)
+(-1)^{p+1}a^0_1\dots a^{p-1}_1.f(a^0_2\ot\dots \ot a^{p-1}_2)a^p+ \]
\[+ \sum_{i=0}^{p-1}(-1)^{i+1}a^0_1\dots a^i_1a^{i+1}_1\dots a^p_1.
f(a^0_2\ot\dots\ot a^i_2.a^{i+1}_2\ot\dots  \ot a^p_2)\]
Using that $f$ is a Hochschild cocycle with values in
$k$, we know that
\[0=
\epsilon(a^0)f(a^1\ot\dots \ot a^p)+
\sum_{i=0}^{p-1}(-1)^{i+1}f(a^0\ot\dots\ot a^i.a^{i+1}\ot\dots  \ot a^p)
+(-1)^{p+1}f(a^0\ot\dots \ot a^{p-1})\epsilon(a^p) \]
So, the summation term in $\partial(F)$ can be replaced using the equality
\[\sum_{i=0}^{p-1}(-1)^{i+1}a^0_1\dots a^i_1a^{i+1}_1\dots a^p_1.
f(a^0_2\ot\dots\ot a^i_2.a^{i+1}_2\ot\dots  \ot a^p_2)=\]
\[= - a^0_1\dots a^i_1a^{i+1}_1\dots a^p_1.\left(
\epsilon(a^0_2)f(a^1_2\ot\dots \ot a^p_2)+
(-1)^{p+1}f(a^0_2\ot\dots \ot a^{p-1}_2)\epsilon(a^p_2) \right)=\]
\[= - \left(a^0.a^1_1\dots a^p_1. f(a^1_2\ot\dots \ot a^p_2)+
(-1)^{p+1}a^0_1\dots a^{p-1}_1.a^pf(a^0_2\ot\dots \ot a^{p-1}_2)\right)
\]
and this finishes the computation of $\partial F$.

Clearly $\epsilon F=f$, so $\epsilon_*$ is a split epimorphism.
To check that $f\mapsto F$ is multiplicative is straightforward:

Let us denote $\wh{f}:=F$, and if $g:A^{\ot q}\to k$, $\wh{g}:A^{\ot q}\to
A$ the cocycle corresponding to $g$.
\[\begin{array}{rcl}
\wh{f\smile g}(a^1\ot\dots \ot a^{p+q})&=&
a^1_1\dots  a^{p+q}_1.(f\smile  g)(a^1_2\ot\dots \ot a^{p+q}_2)\\
&=&
a^1_1\dots  a^{p+q}_1.f(a^1_2\ot\dots \ot a^p_2)g(a^{p+1}_2\ot
\dots\ot a^{p+q}_2)\\
&=&(\wh{f}\smile\wh{g})(a^1\ot\dots \ot a^{p+q})
\end{array}\]
\end{proof}
\section{Brace operations}
In this section we prove our main theorem, stating that
the map $H^{\bullet}(A,k)\to H^{\bullet}(A,A)$ is
``compatible'' with the brace operations, and as a consequence
with the Gerstenhaber bracket.

\begin{teo}\label{teo3}
The image $H^{\bullet}(A,k)\to H^{\bullet}(A,A)$ is stable 
under the brace operation.
Moreover, if $\wh{f}$ and $\wh{g}$ are the images in $H^{\bullet}(A,A)$ 
of $f$ and $g$ in $H^{\bullet}(A,k)$, then 
$\wh{f}\circ_i\wh{g}=\wh{f\circ_i\wh{g}}$
\end{teo}
\begin{proof} Let us recall the definition of the
brace operations (see \cite{G}).
If $F:A^{\ot p}\to H$ and $G:A^{\ot q}\to A$ and $1\leq i\leq p$, then
$F\circ_iG:A^{\ot p+q-1}\to A$ is defined by
\[(F\circ_iG)(a^1\ot \dots\ot a^{i}\ot b^1\ot\dots\ot b^{q}
\ot a^{i+1}\ot\dots\ot a^{p})=
F(a^1\ot\dots\ot a^{i}\ot G(b^1\ot\dots\ot b^{q})
\ot a^{i+1}\ot\dots\ot a^{p})\]
Asume now that $f:A^{\ot p}\to k$, $g:A^{\ot q}\to k$ and
$F=\wh{f}$ and $G=\wh{g}$, namely
\[F(a^1\ot \dots\ot a^{p})=
a^{1}_1\dots a^{p}_1.
f(a^{1}_2\ot \dots\ot a^{p}_2)\] 
and similarly for $G$ and $g$.
Then
\[\begin{array}{l}
(F\circ_iG)(a^1\ot \dots\ot a^{i}\ot b^1\ot\dots\ot b^{q}\ot a^{i+1}\ot
\dots\ot a^{p})\\
=F\left(a^1\ot\dots\ot a^{i}\ot G(b^1\ot \dots\ot b^{q})
\ot a^{i+1}\ot \dots\ot a^{p}\right)\\
=F\left(a^1\ot \dots\ot a^{i}\ot b^1_1\dots b^{q}_1.
g(b^1_2\ot \dots\ot b^{q}_2) \ot a^{i+1}\ot\dots\ot a^{p}\right)\\
= a^{1}_1\dots a^{i}_1.b^{1}_1\dots b^{q}_1.
a^{i+1}_1 \dots a^{p}_1.
f\left( a^{1}_2\ot\dots \ot a^{i}_2\ot b^{1}_2\dots 
b^{q}_2.g(b^{1}_3\ot \dots \ot b^{q}_3)\ot
a^{i+1}_2\ot \dots \ot a^{p}_2\right)\\
=\wh{f\circ_i G}(a^1\ot \dots\ot a^{i}\ot b^1\ot \dots\ot b^{q}
\ot a^{i+1}\ot\dots\ot a^{p})
\end{array}\]
\end{proof}
Recall that the brace operations define a ``composition'' operation
$F\circ G=\sum_{i=1}^{p}(-1)^{q(i-1)}F\circ_iG$, where
$F\in H^p(A,A)$ and $G\in H^q(A,A)$. The Gerstenhaber bracket
is defined as the commutator of this composition, so we
have the desired corollary:
\begin{coro}If $A$ is a Hopf algebra, then
$H^{\bullet}(A,k)$ is a Gerstenhaber subalgebra of $H^{\bullet}(A,A)$.
\end{coro}
Consider $H$ a finite dimensional Hopf algebra
and $X=X(H)$ the algebra defined by C. Cibils and M. Rosso 
(see \cite{CR}). We can prove, at least in the finite dimensional
case, the conjecture of \cite{T3} that $H^{\bullet}_{GS}(H,H)$ is
a Gerstenhaber algebra:
\begin{coro}\label{coroimportante}
Let $H$ be a finite dimensional Hopf algebra,
then $H^{\bullet}_{GS}(H,H)$ $(\cong H_{A4}^{\bullet}(H,H)\cong 
\Ext^{\bullet}_X(H,H))$
is a Gerstenhaber algebra.
\end{coro}
\begin{proof} 
The isomorphism
$H^{\bullet}_{GS}(H,H)\cong H_{A4}^{\bullet}(H,H)
\cong \Ext^{\bullet}_X(H,H)$ was
proved in \cite{T2}. 

Let $A$ denote $D(H)$, the Yetter-Drinfeld
double of $H$. One knows that ${}_X$-$\mod \cong {}_A$-$\mod$,
then $\Ext^{\bullet}_X(H,H)\cong \Ext_A^{\bullet}(H^{coH},H^{coH})=
\Ext_A^{\bullet}(k,k)$, and this a Gerstenhaber subalgebra
of $H^{\bullet}(A,A)$.
\end{proof}

\begin{coro}
Let $H$ be a Hopf algebra and assume that $H$ is a Koszul algebra
(i.e. $H$ is graded, $H_0=k$, and $E(E(H))=H$, where $E(\Lambda)
=\Ext^{\bullet}_{\Lambda}(k,k)$, for an augmented algebra $\Lambda$).
Then $E(H)$ is graded commutative.
\end{coro}


\begin{thebibliography}{30}
\bibitem[Be]{B} D. Benson, Representations and cohomology, Vol I. 
Cambridge Studies in Advanced Mathematics. 30.
Cambridge University Press (1998). 

\bibitem[C-R]{CR} Claude Cibils and  Marc Rosso,
{\em Hopf bimodules are modules}. 
J. Pure Appl. Algebra 128, No.3, 225-231 (1998).

\bibitem[Ge]{G} Murray Gerstenhaber,
{\em The cohomology structure of an associative ring}. 
Ann. Math. (2) 78, 267-288 (1963). 

\bibitem[G-S1]{GS1} Murray Gerstenhaber and Samuel D. Schack,
{\em Bialgebra cohomology, deformations, and quantum groups}. 
Proc. Natl. Acad. Sci. USA 87, No.1, 478-481 (1990).

\bibitem[G-S2]{GS2} Murray Gerstenhaber and Samuel D. Schack,
{\em Algebras, bialgebras, quantum groups, and algebraic deformations}. 
Deformation theory and quantum groups with applications to 
mathematical physics, Proc. AMS-IMS-SIAM Jt. Summer
Res. Conf., Amherst/MA (USA) 1990, Contemp. Math. 134, 51-92 (1992).

\bibitem[\c St]{St}Dragos \c Stefan, 
{\em Hochschild cohomology on Hopf-Galois extensions}.
{ J. Pure Appl. Alg.}, {\bf 103}-2 (1995), pp. 221-233.

\bibitem[Ta1]{T1} Rachel Taillefer, 
{\em Cohomology theories of Hopf bimodules and 
cup-product}. 
C. R. Acad. Sci., Paris, S\'er. I, Math. 332, No.3, 189-194 (2001). 

\bibitem[Ta2]{T2} Rachel Taillefer, thesis Universit\'e Montpellier 2
(2001).

\bibitem[Ta3]{T3} Rachel Taillefer, 
{\em  Injective Hopf bimodules, cohomologies of 
infinite dimensional Hopf algebras and graded-commutativity
of the Yoneda product}. \texttt{ArXivMath math.KT/0207154}.

\end{thebibliography}
\end{document}